\documentclass[11pt]{amsart}
\usepackage{amsfonts,amssymb,graphicx}

\input prepictex
\input pictex
\input postpictex

\numberwithin{equation}{section}

\newtheorem{theorem}{Theorem}[section]
\newtheorem{lemma}[theorem]{Lemma}
\newtheorem{proposition}[theorem]{Proposition}

\theoremstyle{definition}

\theoremstyle{remark}

\newcommand{\CC}{\mathbb{C}}
\newcommand{\HH}{\mathbb{H}}

\newcommand{\ZZ}{\mathbb{Z}}

\newcommand{\hh}{\mathfrak{h}}

\newcommand{\la}{\langle}
\newcommand{\ra}{\rangle}
\newcommand{\one}{\mathbf{1}}
\newcommand{\h}{\mathbf{h}}
\newcommand{\ttt}{\mathfrak{t}}
\begin{document}
\title{Jack polynomials and the coinvariant ring of $G(r,p,n)$.}
\author{Stephen Griffeth}
\address{Department of Mathematics \\
University of Minnesota \\
Minneapolis, MN 55455 \\}
\email{griffeth@math.umn.edu}

\begin{abstract}
We study the coinvariant ring of the complex reflection group $G(r,p,n)$ as a module for the corresponding rational Cherednik algebra $\HH$ and its generalized graded affine Hecke subalgebra $\mathcal{H}$.  We construct a basis consisting of non-symmetric Jack polynomials, and using this basis decompose the coinvariant ring into irreducible modules for $\mathcal{H}$.  The basis consists of certain non-symmetric Jack polynomials, whose leading terms are the ``descent monomials'' for $G(r,p,n)$ recently studied by Adin, Brenti, and Roichman and Bagno and Biagoli.  The irreducible $\mathcal{H}$-submodules of the coinvariant ring are their ``colored descent representations''.
\end{abstract}

\maketitle

\section{Introduction.}

The aim of this paper is to understand the coinvariant ring for the complex reflection group $G(r,p,n)$ as a module over the rational Cherednik algebra and its generalized graded affine Hecke subalgebra.  As applications we construct a basis for the coinvariant ring consisting of certain of the non-symmetric Jack polynomials discovered in \cite{DuOp}, and give a new realization of the ``colored descent representations'' studied in \cite{BaBi} as irreducible modules for a generalized graded affine Hecke algebra.  

The classical formulas
\begin{equation*}
\sum_{w \in S_n} t^{l(w)}=\prod_{i=1}^n \frac{1-t^i}{1-t}=\sum_{w \in S_n} t^{\text{maj}(w)}
\end{equation*} where $S_n$ is the group of permutations of $\{1,2,\dots,n\}$, $l(w)$ is the \emph{length} of $w$ and $\text{maj}(w)$ is the \emph{major index} of $w$ may be obtained by computing the Hilbert series for the coinvariant ring of the symmetric group in three ways: the left hand side corresponds to the divided difference basis, the middle is the quotient of the Hilbert series for all polynomials by the Hilbert series for symmetric polynomials, and the right hand side corresponds to the descent basis.  

In \cite{ABR} Adin, Brenti and Roichman constructed ``colored descent bases'' and a corresponding ``flag major index'' for the coinvariant rings of the type $B$ Weyl groups $G(2,1,n)$ (see \cite{ABR}).  As an application they decompose the coinvariant ring into ``colored descent representations''.  Analagous results for the groups $G(r,p,n)$ are contained in the paper \cite{BaBi} of Bagno and Biagoli.  Our main theorem (\ref{main theorem}) constructs a new basis consisting of non-symmetric Jack polynomials by viewing the coinvariant ring as a module for the rational Cherednik algebra, and shows that upon restriction to the generalized graded affine Hecke algebra inside $\HH$, the coinvariant ring decomposes into irreducible submodules corresponding to ``colored descent classes'' of elements of $G(r,p,n)$.  These are the representations studied without the use of Hecke algebras in \cite{BaBi}, and the leading terms of our basis elements are their ``colored descent monomials''.  We suspect that a version of our results holds for Weyl groups, upon replacing the rational Cherednik algebra and Jack polynomials with the double affine Hecke algebra and Macdonald polynomials, and using the exponential coinvariant ring in place of the coinvariant ring.  

We hope our results may eventually shed some light on the seeming intractibility of the corresponding problems for the exceptional complex reflection groups.  Here the missing ingredient seems to be an analog of the generalized graded affine Hecke subalgebra of the rational Cherednik algebra.  

{\bf Acknowledgements.}  Part of this paper is based on a thesis (\cite{Gri}) written at the University of Wisconsin under the direction of Arun Ram.  I am greatly indebted to him for teaching me about affine Hecke algebras and for suggesting the problem that motivated this work.  

\section{Preliminaries and notation} \label{Prelims}

Let $\hh$ be a finite dimensional complex vector space.  A \emph{reflection} is an element $s \in GL(\hh)$ such that $\text{codim}(\text{fix} (s))=1$.  A \emph{complex reflection group} is a finite subgroup $W$ of $GL(\hh)$ that is generated by the set of reflections it contains.  

Let $W$ be a complex reflection group, let $T$ be the set of reflections in $W$, let $\kappa$ be a variable, and let $c_s$ be a set of variables indexed by $s \in T$ such that $c_{w s w^{-1}}=c_s$ for all $s \in T$ and $w \in W$.  Let $F$ be the field of rational functions with complex coefficients in the variables $\kappa$ and $c_s$, and by abuse of notation write $\hh$ and $\hh^*$ for the vector spaces obtained by extension of scalars to $F$.  Later we will specialize $\kappa=0$; in what follows note that $\kappa=0$ does not make any denominators vanish.

We write
\begin{equation}
F W=F \text{-span} \{t_w \ | \ w \in W \} \quad \text{with multiplication} \quad t_v t_w=t_{vw} \quad \hbox{for $v,w \in W$.}
\end{equation}  The \emph{semi-direct product} of the tensor algebra $T(\hh \oplus \hh^*)$ and $F W$ is
\begin{equation}
T(\hh \oplus \hh^*) \rtimes F W=T(\hh \oplus \hh^*) \otimes_F  FW \quad \text{with multiplication} \quad (f \otimes t_v)(g \otimes t_w) =f v.g \otimes t_{vw}
\end{equation}  for $f,g \in T(\hh \oplus \hh^*)$ and $v,w \in W$.  From now on we will drop the tensor signs.  The \emph{rational Cherednik algebra} is 
\begin{equation}
\HH=TV \rtimes F W/I,
\end{equation} where $I$ is the ideal generated by the relations
\begin{equation} \label{fundamental relation}
yx-xy=\kappa \la x,y \ra-\sum_{s \in T} c_s \la \alpha_s,y \ra \la x,\alpha_s^\vee \ra t_s \quad \hbox{for $x \in \hh^*$, $y \in \hh$,}
\end{equation} and
\begin{equation}
xy=yx \quad \hbox{for $x,y \in \hh$ or $x,y \in \hh^*$.}
\end{equation}  The PBW-theorem for $\HH$ (see \cite{Dri}, \cite{EtGi}, and \cite{RaSh}) asserts 
\begin{equation}
\HH \simeq S(\hh) \otimes S(\hh^*) \otimes F W
\end{equation} as a vector space, where $S(\hh)$ and $S(\hh^*)$ are the symmetric algebras of $\hh$ and $\hh^*$.  It can be proved (for a field of any characteristic) by a straightforward adaptation of the standard proof of the PBW theorem for universal enveloping algebras.

Given a $F W$-module $V$, define the \emph{Verma} module $M(V)$ by
\begin{equation}
M(V)=\text{Ind}_{S(\hh) \otimes F W}^\HH V,
\end{equation} where the set of positive degree polynomials $S(\hh)_{>0}$ acts by $0$ on $V$.  The PBW theorem implies that as a complex vector space
\begin{equation}
M(V) \simeq S(\hh^*) \otimes V.
\end{equation}  In particular, when $V=\one$ is the trivial representation of $F W$, we obtain the \emph{polynomial representation} of $\HH$: 
\begin{equation} \label{DunklopsFormula}
M(\one) \simeq S(\hh^*) \quad \text{with} \quad y.f=\kappa \partial_y f-\sum_{s \in  T} c_s \la \alpha_s,y \ra \frac{f-sf}{\alpha_s}
\end{equation} for $y \in \hh$ and $f \in S(\hh^*)$, where $\partial_y$ denotes the partial derivative in the direction $y$.  These are the famous \emph{Dunkl operators}.  From our point of view, the fact that they commute is a consequence of the PBW theorem, though it is possible to prove the commutativity directly (\cite{DuOp}, for instance).

In Lemma \ref{action lemma} and Theorem \ref{main theorem} we will need the following notation: for $\mu \in \ZZ_{\geq 0}^n$ let $w_\mu \in S_n$ be the shortest permutation such that $w_\mu^{-1}.\mu$ is a partition, and let $v_\mu=w_0 w_\mu^{-1}$ be the longest permutation such that $v_\mu.\mu$ is an anti-partition.

\section{The coinvariant ring of a complex reflection group}

There is a useful \emph{Casimir element} ${\mathbf h}$ in the algebra $\HH$ that helps to distinguish between different lowest weight modules.  This element is the analogue for $\HH$ of the \emph{Euler vector field} $\sum_{i=1}^n x_i \frac{\partial}{\partial x_i}$ in the Weyl algebra.  Fix dual bases $x_1,\dots, x_n$ of $\hh^*$ and $y_1,\dots,y_n$ of $\hh$.  It is straightforward to check that the sum
\begin{equation*}
\sum_{i=1}^n x_i y_i \in \HH
\end{equation*} does not depend on the choice of dual bases.  Let
\begin{equation} \label{element h def}
\h=\sum_{i=1}^n x_i y_i + \sum_{s \in T} c_s (1-t_s).
\end{equation}  We have introduced the shift by $\sum c_s$ in order to simply some formulas that occur later on.  A calculation shows
\begin{equation} \label{commutator for h} 
[\h,x]=\kappa x \quad \text{for} \quad x \in \hh^*, \quad [\h,y]=-\kappa y \quad \text{for} \quad y \in \hh, \quad \text{and} \quad \h t_w=t_w \h \quad \text{for} \quad w \in W,
\end{equation} so that if $\kappa=0$, then $\h$ is \emph{central} in $\HH_c$.  For an irreducible $F W$-module $V$, define $c_V$ to be the scalar by which the element $\sum c_s(1-t_s) \in Z(F W)$ acts on $V$.  Then since the reflections $s \in T$ generate $W$, $V$ is the trivial $F W$ module exactly if $c_V=0$.  In the next proposition we use the fact that if $V \in \text{Irr}(\CC W)$ then $M(V)$ has a unique maximal proper graded submodule (even when $\kappa=0$; otherwise the term ``graded'' may be omitted).  The corresponding irreducible quotient is denoted $L(V)$.  

For real reflection groups, the following proposition is a consequence of the results in \cite{Dun}.

\begin{proposition} \label{irred prop}
Suppose that $\kappa=0$ but the other parameters $c_s$ remain generic.  Let $I=S(\hh^*)^W_+ S(\hh^*)$ be the ideal generated by the positive degree $W$-invariant polynomials.  Then the coinvariant ring $S/I$ is an irreducible $\HH$-module.
\end{proposition} 
\begin{proof}
In light of our assumption that $\kappa=0$ and \eqref{DunklopsFormula}, the ideal $I$ is $\HH$-stable and the coinvariant ring is an $\HH$-module. Let $R$ be the (unique) maximal proper graded submodule of $M(\one)$.  We must show that $R \subseteq I$.  It suffices to prove that for all irreducible $\CC W$-modules $V$ and all integers $d$ that $(R^d)^V$ (the $V$-isotypic component of $R^d$) is contained in $I$.  Suppose towards a contradiction that this is false, and choose $d$ minimal so that it fails.  Let $f \in (R^d)^V$ and suppose $f \notin I$.  Note that $d>0$ since $R^0=0$.  Calculate
\begin{equation*}
0=\h.f=\left( \sum_{i=1}^n x_i y_i - \sum_{s \in T} c_s (1-t_s) \right).f=\sum_{i=1}^n x_i y_i.f-c_V f.
\end{equation*}  By minimality of $d$, $\sum_{i=1}^n x_i y_i .f \in I$, so that
\begin{equation}
c_V f \in I.
\end{equation}  Thus $c_V=0$ and our hypothesis implies $V=\one$, contradiction.
\end{proof}    

Our strategy in the remainder of the paper is to construct a basis of the coinvariant ring for the groups $G(r,p,n)$ that is particularly adapted to understanding its structure as an $\HH$-module.

\section{Non-symmetric Jack polynomials}

From now on we consider the case $W=G(r,p,n)$, where $G(r,p,n)$ is the group of $n$ by $n$ monomial matrices whose non-zero entries are $r$th roots of $1$ and so that the product of the non-zero entries is an $r/p$th root of $1$.  We put $\zeta=e^{2 \pi i/r}$ and let $\zeta_i$ be the diagonal matrix with a $\zeta$ in the $i$th position and $1$'s elsewhere.  As usual, $s_{ij}$ is the transposition matrix interchanging the $i$th and $j$th coordinates.

 When $n \geq 3$ the equations
\begin{equation}
(\zeta_i^l \zeta_k^{-l}) s_{ij} (\zeta_i^l \zeta_k^{-l})^{-1}=\zeta_i^l s_{ij} \zeta_i^{-l} \quad \text{and} \quad s_{1i} \zeta_i^l s_{1i}^{-1}=\zeta_1^l,
\end{equation} for $1 \leq i<j \leq n$, $k \neq i,j$, and $0 \leq l \leq r-1$, show that there are $r/p$ conjugacy classes of reflections in $G(r,p,n)$:
\begin{itemize}
\item[(a)]  The reflections of order two:
\begin{equation}
\zeta_i^l s_{ij} \zeta_i^{-l}, \quad \text{for} \quad 1 \leq i < j \leq n, \quad 0 \leq l \leq r-1, 
\end{equation} and
\item[(b)]  the remaining $r/p-1$ classes, consisting of diagonal matrices
\begin{equation}
\zeta_i^{pl}, \quad \text{for} \quad 1 \leq i \leq n, \quad 1 \leq l \leq r/p-1,
\end{equation} where $\zeta_i^{pl}$ and $\zeta_j^{pk}$ are conjugate if and only if $k=l$.
\end{itemize}  Despite the fact that there are more conjugacy classes than described above when $n=2$ and $p>1$, the results in this paper go through without change.  

Let
\begin{equation*}
y_i=(0,\dots,1,\dots,0)^t \quad \text{and} \quad x_i=(0,\dots,1,\dots,0)
\end{equation*} have $1$'s in the $i$th position and $0$'s elsewhere, so that $y_1,\dots,y_n$ is the standard basis of $\hh=\CC^n$ and $x_1,\dots,x_n$ is the dual basis in $\hh^*$.  

By translating the definition given in Section \ref{Prelims}, the \emph{rational Cherednik algebra} $\HH$ for $G(r,1,n)$ with parameters $\kappa,c_0,c_1,\dots,c_{r-1}$ is the $F$-algebra generated by $F[x_1,\dots,x_n]$, $F[y_1,\dots,y_n]$, and $F G(r,1,n)$ with relations 
\begin{equation*}
t_w x=(wx) t_w \quad \text{and} \quad t_w y=(wy) t_w,
\end{equation*} for $w,v \in W$, $x \in \hh^*$, and $y \in \hh$, 
\begin{equation} \label{cf2}
y_i x_j=x_j y_i+c_0 \sum_{l=0}^{r-1} \zeta^{-l} t_{\zeta_i^l s_{ij} \zeta_i^{-l}},
\end{equation} for $1\leq i \neq j \leq n$, and 
\begin{equation} \label{cf3}
y_i x_i=x_i y_i+\kappa-\sum_{l=1}^{r-1} c_l (1-\zeta^{-l}) t_{\zeta_i^l} -c_0 \sum_{j \neq i} \sum_{l=0}^{r-1} t_{\zeta_i^l s_{ij} \zeta_i^{-l}},
\end{equation} for $1 \leq i \leq n$.  

In the above definition, suppose $c_l=0$ for $l$ not divisible by $p$.  The \emph{rational Cherednik algebra} $\HH$ for $W=G(r,p,n)$ with parameters $\kappa,c_0,c_p,\dots,c_{r-p}$ is the subalgebra of the rational Cherednik algebra for $G(r,1,n)$ generated by $G(r,p,n)$, $F[x_1,\dots,x_n]$, and $F[y_1,\dots,y_n]$.  

In section 3 of \cite{DuOp} a very useful commutative subalgebra of $\HH$ is defined: let

\begin{equation}
z_i=y_i x_i+c_0 \phi_i \quad \text{where} \quad \phi_i=\sum_{1 \leq j <i} \sum_{l=0}^{r-1} t_{\zeta_i^l s_{ij} \zeta_i^{-l}}.
\end{equation}  Let $\ttt$ be the commutative subalgebra of $\HH$ generated by $z_1,\dots,z_n$, $t_{\zeta_1}^p,\dots,t_{\zeta_n}^p$ and $t_{\zeta_1^{-1} \zeta_{2}},\dots,t_{\zeta_{n-1}^{-1} \zeta_n}$.  If $\alpha: \ttt \rightarrow F$ is an $F$-algebra homomorphism and $M$ is an $\HH$-module, define the $\alpha$-\emph{weight space} $M_\alpha$ by
\begin{equation}
M_\alpha=\{m \in M \ | \ \hbox{there is $q \in \ZZ_{>0}$ such that }(f-\alpha(f))^q.m=0 \quad \hbox{for all $f \in \ttt$.} \}
\end{equation}   If $v \in M$, we say that $v$ has $\ttt$-weight $(\alpha_1,\dots,\alpha_n,\zeta^{\beta_1},\dots,\zeta^{\beta_n})$ if
\begin{equation} \label{weight conditions}
z_i.v=\alpha_i v, \quad t_{\zeta_i^p}.v=\zeta^{p \beta_i}, \quad \text{and} \quad t_{\zeta_i^{-1} \zeta_{i+1}}.v=\zeta^{\beta_{i+1}-\beta_i} v.
\end{equation}  If $v \neq 0$ then the sequence $(\alpha_1,\dots,\alpha_n,\zeta^{\beta_1},\dots,\zeta^{\beta_n})$ is determined by \eqref{weight conditions} up to simultaneously multiplying $\zeta^{\beta_1},\dots,\zeta^{\beta_n}$ by a power of $\zeta^{r/p}$.  Then by the case $\lambda=(n)$ of Theorem 5.1 in \cite{Gri2} there is a unique basis $f_\mu$ of $S(\hh^*)=\CC[x_1,\dots,x_n]$ such that $f_\mu$ is a $\ttt$-eigenvector and
\begin{equation}
f_\mu=x^\mu+\text{lower terms},
\end{equation} where the lower terms are with respect to a certain partial order on $\ZZ_{\geq 0}^n$ extending dominance order on partitions.  In fact, $f_\mu$ is essentially a non-symmetric Jack polynomial; see Proposition 3.14 of \cite{DuOp} and the material preceding it.  

The \emph{generalized degenerate affine Hecke algebra} is the subalgebra $\mathcal{H}$ of $\HH$ generated by $\CC W$ and $\ttt$.  It was first constructed in \cite{RaSh}, section 5, and in \cite{Dez} it was observed that by the results of \cite{DuOp} it is a subalgebra of $\HH$.  For $W=G(1,1,n)$, it is the usual graded affine Hecke algebra of the symmetric group.

Some formulas are simpler when written in terms of the following parameters:
\begin{equation} \label{ds def}
d_j=\sum_{l=1}^{r/p-1} c_{lp} \zeta^{lpj}, \quad \hbox{for $j \in \ZZ/ r \ZZ$. }
\end{equation}  To efficiently describe the $\HH$-action on the basis $f_\mu$, we introduce the following operators:
\begin{equation} \label{sigma def}
\sigma_i=t_{s_i}+\frac{c_0}{z_i-z_{i+1}} \pi_i \quad \hbox{for $1 \leq i \leq n-1$} \quad \hbox{where $\pi_i=\sum_{l=0}^{r-1} t_{\zeta_i \zeta_{i+1}^{-1}}^l$.}
\end{equation}  The operator $\sigma_i$ is well-defined on those $\ttt$-weights spaces $M_\alpha$ on which $z_i-z_{i+1}$ is invertible or $\pi_i$ acts as $0$.  We also define the intertwining operators $\Phi$ and $\Psi$ by
\begin{equation} \label{phipsi defs}
\Phi=x_n t_{s_{n-1} \cdots s_1} \quad \text{and} \quad \Psi=y_1 t_{s_1 \cdots s_{n-1}}.
\end{equation}  The intertwiner $\Phi$ was discovered by Knop and Sahi (\cite{KnSa}).  

To calculate the action of the intertwiners on $f_\mu$, we record the $\ttt$-eigenvalue of $f_\mu$.  By Theorem 5.1 of \cite{Gri}, it is given by
\begin{equation} \label{eig1}
z_i.f_\mu=\left((\mu_i+1) \kappa -(d_0-d_{-\mu_i-1})-r(v_\mu(i)-1) c_0 \right) f_\mu,
\end{equation}
\begin{equation} \label{eig2}
t_{\zeta_i^p}.f_\mu=\zeta^{-p \mu_i} f_\mu \quad \text{and} \quad t_{\zeta_i \zeta_{i+1}^{-1}}.f_\mu=\zeta^{\mu_{i+1}-\mu_i} f_\mu.
\end{equation}  We also need the following operators on multi-indices $\mu \in \ZZ_{\geq 0}$:
\begin{equation*}
\phi.(\mu_1,\mu_2,\dots,\mu_n)=(\mu_2,\mu_3,\dots,\mu_n,\mu_1+1) \quad \text{and} \quad \psi.(\mu_1,\mu_2,\dots,\mu_n)=(\mu_n-1,\mu_1,\mu_2,\dots,\mu_{n-1}).
\end{equation*}

The following lemma gives the action of the intertwiners on the basis $f_\mu$ and is a special case of Lemma 5.2 of \cite{Gri2}.

\begin{lemma} \label{action lemma}
Let $\mu \in \ZZ_{\geq 0}^n$.
\begin{enumerate}
\item[(a)] If $\mu_i<\mu_{i+1}$ or $\mu_i-\mu_{i+1} \neq 0$ mod $r$ then
\begin{equation*}
\sigma_i.f_{\mu}=f_{s_i.\mu}.
\end{equation*}
\item[(b)]  If $\mu_i>\mu_{i+1}$ and $\mu_i-\mu_{i+1} =0$ mod $r$ then
\begin{equation*}
\sigma_i.f_{\mu}=\frac{(\delta-rc_0)(\delta+rc_0)}{\delta^2} f_{s_i.\mu},
\end{equation*} where
\begin{equation*}
\delta=\kappa(\mu_{i}-\mu_{i+1})-c_0 r (v_\mu(i)-v_\mu(i+1)).
\end{equation*}
\item[(c)]  For all $\mu \in \ZZ_{\geq 0}^n$,
\begin{equation*}
\Phi.f_{\mu}=f_{\phi.\mu}.
\end{equation*}
\item[(d)]  For all $\mu \in \ZZ_{\geq 0}^n$,
\begin{equation*}
\Psi.f_{\mu}=\begin{cases}  \left(\kappa \mu_n-(d_0-d_{-\mu_n})-c_0 r (v_\mu(n)-1)  \right) f_{\psi.\mu}\quad &\hbox{if $\mu_n>0$,} \\ 0 \quad &\hbox{if $\mu_n=0$.} \end{cases}
\end{equation*}
\end{enumerate}
\end{lemma}

As a consequence of this lemma, the polynomials $f_\mu$ are well-defined at $\kappa=0$: they can be recursively constructed by using $\Phi$ and the $\sigma_i$'s, which are well-defined on $M(\one)$ when $\kappa=0$. 

\section{The coinvariant ring of $G(r,p,n)$.}

In this section we assume $\kappa=0$ and that $W=G(r,p,n)$.

We will now obtain an eigenbasis for the coinvariant ring $S/I$ for $G(r,p,n)$ indexed by a certain subset of $G(r,1,n)$.  First we need some definitions.  We write elements of $G(r,1,n)$ as ``colored permutations'':
\begin{equation}
v=\left[\zeta^{k_1} w(1), \zeta^{k_2} w(2), \dots, \zeta^{k_n} w(n) \right]=w \zeta^{k_1} \cdots \zeta^{k_n}, 
\end{equation} with $w \in S_n$, $0 \leq k_1,\dots,k_n \leq r-1$.  A \emph{descent} of $v$ is an integer $1 \leq i \leq n-1$ such that
\begin{equation}
k_i < k_{i+1} \qquad \text{or} \qquad \hbox{$k_i=k_{i+1}$ and $w(i)>w(i+1)$.}
\end{equation}  The \emph{Steinberg weight} for $v$ is $\lambda_v=(d_1(v),\dots,d_n(v))$, where
\begin{equation}
d_i(v)=r \left|\{j \geq w^{-1}(i) \ | \ \hbox{$j$ is a descent of $v$} \}\right|+k_{w^{-1}(i)},
\end{equation} and $w^{-1}(i)$ is the position of $i$ in the sequence $[w(1),\dots,w(n)]$.  The \emph{colored descent class} of $v \in G(r,1,n)$ is the pair $\text{des}(v)=(d(v),(k_1,k_2,\dots,k_n))$, where $d(v)$ is the set of positions $1 \leq i \leq n-1$ in which $v$ has a descent.  We write $D_p$ for the set of all descent classes of elements of $G(r,1,n)$ satisfying $k_n \leq r/p-1$.

The \emph{colored descent monomial} corresponding to $v$ is
\begin{equation}
x^{\lambda_v}=x_1^{d_1(v)} x_2^{d_2(v)} \dots x_n^{d_n(v)},
\end{equation}  If $w_\mu$ is the shortest permutation such that $w_\mu^{-1}.\mu$ is a partition, then we note that 
\begin{equation} \label{descent facts}
w=w_{\lambda_v}, \quad k_{w^{-1}(i)}=d_i(v) \ \hbox{mod  $r$}, \quad \text{and} \quad s_i.\lambda_v=\lambda_{s_i v} \ \hbox{if $\text{des}(v)=\text{des}(s_iv)$.}
\end{equation}  We will also need the formulas
\begin{equation} \label{eig3}
z_i.f_\mu=(d_{-\mu_i-1}-d_0-r(v_\mu(i)-1) c_0) f_\mu, \quad t_{\zeta_i^p}.f_\mu=\zeta^{-p \mu_i} f_\mu, \quad \text{and} \ t_{\zeta_i \zeta_{i+1}^{-1}}.f_\mu=\zeta^{-(\mu_i-\mu_{i+1})}f_\mu
\end{equation} obtained by specializing \eqref{eig1} and \eqref{eig2} to $\kappa=0$.

Our proof our \ref{main theorem} requires the following combinatorial lemma:
\begin{lemma}\label{chains lemma}
Let $v,v' \in G(r,1,n)_p$ with $\text{des}(v)=\text{des}(v')$.  Then there exists a sequence $s_{i_1},\dots,s_{i_q}$ of simple transpositions so that $\text{des}(s_{i_j} \cdots s_{i_1} v)=\text{des}(s_{i_{j-1}} \dots s_{i_1} v)$ for $1 \leq j \leq q$ and $v'=s_{i_q} \cdots s_{i_1} v$.  
\end{lemma}
\begin{proof}
Each descent class contains a unique $v=w \zeta_1^{k_1} \cdots \zeta_n^{k_n}$ so that if $i$ is a descent of $v$ with $k_i=k_{i+1}$ then $w(i)=w(i+1)+1$ and otherwise $w(i)<w(i+1)$.  One checks that such a $v$ is the unique element of its descent class with $w$ of minimum length, and it is straightforward to check that for any other $v'=w' \zeta_1^{k_1} \cdots \zeta_n^{k_n}$ in the same descent class, there is a simple reflection $s_i$ with $l(s_i w')<l(w')$ and $\text{des}(s_i v')=\text{des}(v')$.  
\end{proof}

The $x^{\lambda_v}$ generalize the descent monomials from \cite{GaSt} and \cite{Gar}, and recently appeared in \cite{BaBi}.  The following theorem shows that they are the leading terms of a basis for the coinvariant ring consisting of certain $\kappa=0$ specializations of non-symmetric Jack polynomials.  It strengthens Theorem 8.8 of \cite{Gri3}.

\begin{theorem} \label{main theorem}
Suppose $\kappa=0$ and $c_s$ are generic.  Then $L(\one)$ is the coinvariant ring for $G(r,p,n)$ and has basis $\{f_{\lambda_v} \ | \ v \in G(r,1,n)_p \}$, where
\begin{equation*} \label{Gr1np def}
G(r,1,n)_p=\left\{ \left[\zeta^{k_1} w(1),\dots,\zeta^{k_n} w(n) \right] \in G(r,1,n) \ | \ 0 \leq k_n \leq r/p-1 \right\}.
\end{equation*}  Furthermore, as a module for the generalized graded affine Hecke algebra $\mathcal{H}$, $L(\one)$ decomposes into irreducibles as
\begin{equation*}
L(\one)=\bigoplus_{d \in D_p} \CC \{ f_{\lambda_v} \ | \ v \in G(r,1,n)_p \ \text{and} \ \text{des}(v)=d \},
\end{equation*} no two of which are isomorphic.
\end{theorem}
\begin{proof}
By Proposition \ref{irred prop} we already know that $L(\one)$ is the coinvariant ring.   The algebra $\HH$ is generated by $\Phi$, $\Psi$, and $\CC W$.  Thus if the span of the $f_\mu$'s such that $\mu$ is \emph{not} a Steinberg weight $\lambda_v$ is stable under the intertwining operators from Lemma \ref{action lemma}, dimension considerations prove the first assertion of the Theorem.  This is checked in a straighforward way using  Lemma \ref{action lemma}, \eqref{descent facts}, and the equivalences
\begin{equation} \label{equiv1}
\sigma_i^2.f_{\lambda_v}=0 \quad \iff \quad k_i=k_{i+1} \ \text{and} \ w^{-1}(i)=w^{-1}(i+1) \pm 1 \quad \iff \quad \text{des}(v) \neq \text{des}(s_i v)
\end{equation} and
\begin{equation}
\Psi \Phi f_{\lambda_v}=0 \quad \iff \quad w^{-1}(1)=n \ \text{and} \ k_n=r/p-1.
\end{equation}

 The formulas \eqref{eig3} show that the $\ttt$-eigenspaces on the span of $f_{\lambda_v}$ for $v \in G(r,1,n)_p$ are all one-dimensional.  Therefore each irreducible $\mathcal{H}$ submodule of $L(\one)$ is spanned by the $f_\mu$'s it contains.  Since $\mathcal{H}$ is generated by $\ttt$ and $s_1,\dots,s_n$, the $\CC$-span of a collection of $f_{\lambda_v}$'s is an $\mathcal{H}$-module exactly if it is stable under $\sigma_1,\dots,\sigma_n$, and is irreducible exactly if any two $f_{\lambda_v}$'s can be connected by a sequence of invertible intertwiners.  On the other hand, Lemma \ref{action lemma}, Lemma \ref{chains lemma}, \eqref{descent facts}, and \eqref{equiv1} can be combined once more to see that if $\text{des}(v)=\text{des}(s_i v)$ then $\sigma_i.f_{\lambda_v}=c f_{\lambda_{s_i v}}$ with $c \neq 0$, and if $v,v' \in G(r,1,n)_p$ then there is a sequence of invertible intertwiners $\sigma_i$ connecting $f_{\lambda_v}$ and $f_{\lambda_{v'}}$ exactly if $\text{des}(v)=\text{des}(v')$.  The second assertion of the Theorem follows from this.  Finally, the summands are pairwise non-isomorphic because their $\ttt$-spectra are different.
\end{proof}  

The above theorem may help explain why the major index is difficult to define directly on the group $G(r,p,n)$: it is the subset $G(r,1,n)_p$ of $G(r,1,n)$ that naturally indexes the basis of Jack polynomials, not the group $G(r,p,n)$.

\bibliographystyle{amsplain}
\def\cprime{$'$} \def\cprime{$'$}

\end{document}